\magnification=1200
\font\bb=msbm10
\def\square{\hbox{\rlap{$\sqcap$}$\sqcup$}}
\def\met{\hbox{Riem}}
\def\ric{\hbox{Riem}^{Ric+}}
\def\scal{\hbox{Riem}^{scal+}}
\centerline{\bf ON THE MODULI SPACE OF POSITIVE RICCI}
\centerline{\bf CURVATURE METRICS ON HOMOTOPY SPHERES}
\bigskip\bigskip
\centerline{\bf David J. Wraith}
\bigskip\bigskip
\itemitem{} {\bf Abstract:} {\it We show that the moduli space of Ricci positive metrics on a certain}
\itemitem{} {\it family of homotopy spheres has infinitely many components.} 
\footnote{}{2000 {\it Mathematics Subject Classification:} 53C20.}
\footnote{}{{\it Keywords:} positive Ricci curvature, moduli space, homotopy spheres}
\bigskip\bigskip
\centerline{\bf \S1 Introduction}
\bigskip

A fundamental problem in Riemannian geometry is to understand metrics
which display certain types of curvature characteristics. Of particular
imporance among curvature characteristics are the various forms of positive
curvature, such as positive sectional curvature, positive Ricci curvature and
positive scalar curvature. Existence issues for positive scalar curvature metrics
are reasonably well understood. The situation for positive Ricci and positive
sectional curvature metrics is somewhat less clear.

For the most part the existence question, that is, whether or not a particular
manifold or family of manifolds admits metrics of a given type, has been a primary 
focus of research. However, there is an equally intruiging secondary question. If
a manifold or family of manifolds admit a given type of metric, how are such
metrics distributed among all possible Riemannian metrics on these objects? For
example are they rare or common?

A particularly natural question concerns connectedness: given two metrics of
a given type on a particular manifold, is it possible to smoothly deform
one metric into the other {\it through metrics of the same type.} If the
answer to this question is always yes, this shows that the set of metrics
having the given type form a connected subset of the space of all metrics
(at least if the latter is equipped with the smooth topology). Thus it
makes sense to ask about the connectedness of the subset of `nice' metrics
in whatever context one is working. If this set is not connected, how many
connected components are there? Of course, one could also seek to uncover more subtle aspects
of the topology, such as the fundamental or higher homotopy groups. For some recent work in this
direction, see [BHSW].

In this paper we will be interested in the subset of metrics with positive scalar and
positive Ricci curvature. If a given manifold admits a positive Ricci curvature
metric, then this metric also has positive scalar curvature. It makes sense then
to ask about the {\it relative} distribution of the positive scalar versus positive
Ricci curvature metrics. For example, does every connected component of positive
scalar curvature metrics (assuming there is more than one) contain a positive Ricci
curvature metric?

Instead of working with spaces of metrics, we might instead choose to work with
{\it moduli spaces} of metrics. For any smooth manifold $M,$ we can consider the
group of all diffeomorphisms of $M,$ $\hbox{Diff}(M).$ Let $\met (M)$ denote
the space of all Riemannian metrics on $M$ with the smooth topology. Then $\hbox{Diff}(M)$
acts on $\met (M)$ by pulling back metrics. The quotient of $\met (M)$ by this
action is a moduli space of metrics. Denoting the space of Ricci positive metrics on $M$
by $\ric (M)$ and the space of positive scalar metrics by $\scal (M),$ we can also
form the moduli spaces of positive Ricci and positive scalar curvature metrics in the
same way: $\ric (M)/\hbox{Diff}(M),$ $\scal (M)/\hbox{Diff}(M).$ Note that a
lower bound on the number of components in a moduli space is also a lower bound on the
number of components for the corresponding space of metrics.

Possibly the earliest result on the connectedness of spaces of metrics displaying some form of positive
curvature was due to Hitchin [H]. He proved that if a closed spin manifold admits a positive scalar curvature
metric and has dimension a multiple of eight, then its space of positive scalar curvature metrics has more than
one path component and has non-trivial fundamental group.

A few years later, R. Carr [C] proved that the space of positive scalar
curvature metrics on the spheres $S^{4k-1}$ has {\it infinitely} many
path-components. This result demonstrated for the first time how complex the distrubution
of positive scalar curvature metrics can be. Of course spheres also admit metrics with positive
Ricci and positive sectional curvature. Thus it is natural to ask whether similar results hold
if scalar curvature is replaced by Ricci or sectional curvature in Carr's result.

The objects we will study in this paper are homotopy spheres. These are smooth manifolds
having the same homotopy type as some (standard) sphere $S^n.$ At least in dimensions other
than four, homotopy spheres are either standard or {\it exotic}, that is, homeomorphic but
not diffeomorphic to the standard sphere. The set of all homotopy spheres has a particularly
nice and very large subfamily: the homotopy spheres which bound parallelisable manifolds. It was
shown in [W1] that these manifolds all admit Ricci positive metrics. (See also [BGN]. For a general
reference about the construction and curvature of homotopy spheres, see [JW].)
Our main result is as follows:

\proclaim Theorem A. Let $\Sigma^{4k-1}$ be any homotopy sphere bounding a parallelisable manifold,
where $k > 1$. Let $\met (\Sigma)$ denote the space of all smooth Riemannian metrics on $\Sigma$ equipped
with the $C^{\infty}$ topology, and denote by $\ric (\Sigma)$ the subset consisting of all Ricci
positive metrics. Give $\ric (\Sigma)$ the induced topology from $\met (\Sigma)$. Then the moduli space 
$\ric (\Sigma)/\hbox{Diff}(\Sigma)$ has infinitely many components.
\par

As an immediate corollary we see that $\ric (\Sigma)$ must also have infinitely many components.
Note that it does {\it not} follow from this that either $\scal (\Sigma)$ or the moduli space $\scal (\Sigma)/\hbox{Diff}(\Sigma)$ must have infinitely
many components, as to deduce this would require further information about the distribution of positive Ricci curvature
metrics among the positive scalar curvature metrics. In fact, both of these positive scalar curvature spaces {\it do}
have infinitely many components.
For example it is not difficult to see (though not pointed out in [C]) that Carr's arguments
apply equally well - with the same conclusion - to any exotic sphere in these dimensions which bounds a parallelisable manifold.

In [KS], M. Kreck and S. Stolz strengthened Carr's result considerably. They showed that
for any closed spin manifold $M^{4k-1}$ ($k > 1$) which admits a positive scalar
curvature metric and for which $H^1(M;\hbox{\bb Z}/2)=0$, the moduli space of positive scalar
curvature metrics has infinitely many components. They then used this result to prove
that the moduli space of Ricci positive metrics on certain homogeneous 7-manifolds with $SU(3) \times SU(2) \times U(1)$ symmetry has infinitely many components. To the best of our knowledge, apart from Theorem A above, these are the only examples of Ricci positive manifolds in the literature for which
the moduli space of Ricci positive metrics has been shown to have infintely many components. As the Kreck-Stolz examples all occur in dimension 7,
Theorem A also demonstrates for the first time that the following corollary is true:

\proclaim Corollary B. There are manifolds in infinitely many dimensions for which the moduli space of Ricci positive metrics has 
infintely many components. \par 

Other results about spaces or moduli spaces of metrics displying some kind of positive curvature are few and far between in the literature.
Shortly after [KS], Botvinnik and Gilkey [BG] showed that in some cases, the space of positive scalar curvature metrics for odd dimensional spin manifolds with non-trivial finite fundamental group can also have infinitely many path components.
More recent results due to Kapovitch, Petrunin and Tuschmann [KPT] show that there are examples of non-compact manifolds in every dimension $\ge 22$ for which the moduli space of non-negatively curved metrics has infinitely many components. Let us also mention that in [KS], examples are given of
Wallach spaces for which the moduli space of non-negatively curved metrics is not connected. It is observed
however in [KPT] that the $SU(3) \times SU(2) \times U(1)$-homogeneous spaces studied in [KS] actually have moduli 
spaces of non-negatively curved metrics with infinitely many components. Most recently, in [BHSW] it is shown
that for many degrees in a stable range, the homotopy groups of the moduli space of positive scalar curvature metrics on 
$S^n$ and other manifolds are non-trivial.
\smallskip

It is interesting to compare and contrast the techniques of Carr and examples of Kreck-Stolz with those of this paper.

Carr regards the homotopy spheres in question as the boundaries of `plumbed' manifolds. Plumbing is a construction whereby disk
bundles are glued together to create a new manifold with boundary. The basic idea is as follows.
Consider two disc bundles, and for each bundle choose
a locally trivial neighbourhood over a disc in the base. We glue the bundles by identifying these neighbourhoods. 
To do this we use a diffeomorphism which identifies the base disc of one neighbourhood with the fibre disc of the other, 
and vice versa. (Of course, this can only work if the dimensions match.) 
The resulting object is the so-called plumbing of the disc bundles, and can be made differentiable
by simply straightening out the angles. By plumbing certain sequences of disc bundles over spheres together, one can construct
every homotopy sphere which bounds a parallelisable manifold (and some others too). The key point is that a given homotopy
sphere of this type can be realised as the boundary of infinitely many different plumbed manifolds.
(For the basics of plumbing and applications to homotopy sphere construction see [Br; V.2] or [W1].)

Carr proves that for any regular neighbourhood of the union
of base spheres in the plumbing construction, there is a metric of positive scalar curvature which is a product near the boundary,
and thus the boundary itself has positive scalar curvature. If two boundary metrics arising in this way from different plumbing
descriptions of the same homotopy sphere belong to the same path component of positive scalar curvature metrics, it is not
difficult to see that the two plumbed manifolds can be glued toether to give a smooth positive scalar curvature manifold.
However, the existence of such a metric on this new object can be ruled out using simple index theory. Thus the 
metrics arising in this way from different plumbing descriptions must all lie in different components of the space
of positive scalar curvature metrics.

The same approach fails to work for positive Ricci curvature as the boundary of the regular neighbourhood will almost certainly
not have positive Ricci curvature. Instead, one needs to {\it begin} with a positive Ricci curvature metric on the boundary, and
hope to extend it inwards over the whole plumbed manifold in (at least) a positive scalar fashion so that a neighbourhood
of the boundary is a product. If one could do this, then the rest of Carr's argument would still apply, showing that the space
of positive Ricci curvature metrics has infinitely many compoments. However, such a construction must necessarily depend on the 
precise form of the Ricci positive metric. A significant part of this paper is dedicated to performing this kind of metric
extension. 

Of course, our main result concerns moduli spaces of Ricci positive metrics, as opposed to just the set of
such metrics. To obtain results about moduli spaces we use the Kreck-Stolz techniques, for which the existence of our metric extension
is in fact a pre-requisite. Kreck and Stolz define a {\bb Q}-valued quantity $s$ for $(4n-1)$-dimensional spin manifolds with positive
scalar curvature and vanishing real Pontrjagin classes, with the property that $|s|$ is an invariant of the moduli space of positive
scalar curvature metrics. The computation of this invarant is difficult in general. However, for the examples we are interested in,
we exploit the fact that our homotopy spheres bound parallelisable manifold to simplify the calculation. In a similar fashion to
the Carr approach, a pre-requisite for
this method to work is that the boundary metric should be extendable over the interior of the plumbed manifold 
in (at least) a positive scalar fashion, in such a way that neighbourhood of the boundary is a product.
Again, different plumbing descriptions of the same homotopy sphere lead to metrics with different $s$-invariants. Thus
the metrics must all lie in different components of the moduli space of positive scalar curvature metrics, and hence in different components
of the moduli space of positive Ricci curvature metrics.

The is a problem, however. The (reasonably explicit) Ricci positive metrics on homotopy spheres constructed in [W1] and [W2] do not
easily extend over the bounding manifold in the desired manner. Our solution is to deform these boundary metrics within
positive scalar curvature to a metric which admits the required extension. 

The metric extensions are constructed in pieces, corresponding to the various disc bundles being plumbed together.
On the boundary, each of these plumbings is equivalent to a surgery, and thus each plumbing can be viewed as the
addition of a handlebody. The metrics on our homotopy spheres are constructed directly on the boundary of this construction,
surgery by surgery. Extending the metric in the appropriate way over the corresponding handlebody is a delicate issue, which
we will explore in detail.
\smallskip

This paper is organised as follows. In \S2 we discuss the construction of homotopy spheres as boundaries
of plumbed manifolds and introduce the s-invariant of Kreck and Stolz. 
In \S3 we look at surgery on Ricci
positive manifolds. In \S4 we show how to deform a positive Ricci curvature metric on an exotic sphere through positive
scalar curvature metrics into a form suitable for extension across the corresponding bounding manifold. This extension is addressed in \S5. 
In \S6 we consider paths of positive scalar curvature metrics, concluding with the proof of Theorem A.

The author would like to thank C. B\"ohm, W. Tuschmann and B. Botvinnik for helpful comments during
the preparation of this article, and especially S. Bechtluft-Sachs for his careful reading and criticism
of earlier versions. Finally, it is a pleasure to thank Mark Walsh for his advice on drawing the figures. 

\bigskip\bigskip
\centerline{\bf \S2 Homotopy Spheres, Plumbing and the s-invariant}
\bigskip
The diffeomorphism classes of homotopy spheres bounding parallelisable manifolds
of dimension $n$ form an abelian group under the connected sum operation. This group is denoted
$bP_n$. It was shown in [KM] that for $n$ odd, $bP_n=0$, for $n=4k+2$, $bP_n$ is either 0
or $\hbox{\bb Z}_2$, and for $n=4k$, $bP_n$ is cyclic.

Let $\Sigma \in bP_{4k}$ for some $k>1$. It was shown in [W1] that $\Sigma$ arises as the boundary
of a plumbed manifold: in fact $\Sigma$ arises as the boundary of a plumbed manifold in
infinitely many different ways. For this plumbed manifold, the disc-bundles involved
are all $D^{2k}$-bundles over $S^{2k}.$
To understand the details of the construction, we can represent the plumbing by a schematic diagram in the following
way. For each bundle we draw a dot, which should be labelled with the appropriate element of
the group $\pi_{2k} BSO(2k)$ which classifies such bundles. Each time we plumb two of the bundles together, join the appropriate dots with a line.
In this way we  construct a graph. Consider the $E_8$ graph as a plumbing graph, where each dot
represent the tangent disc bundle over $S^{2k}.$  
\bigskip
\input rlepsf.tex
\centerline{
\relabelbox \epsfxsize 1.5 truein
\epsfbox{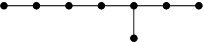}
\endrelabelbox} 
\medskip
\centerline{\it Figure 1: the $E_8$ graph.}
\medskip
\noindent The boundary of the manifold resulting from this plumbing can be shown to represent a generator of $bP_{4k}$
(see [Br]). We obtain the $n$th multiple of this generator by starting with the $E_8$ graph and successivly attaching $n-1$ copies of the following
$\lq$building block' to the right-hand end.
\bigskip
\centerline{
\relabelbox \epsfxsize 1.5 truein
\epsfbox{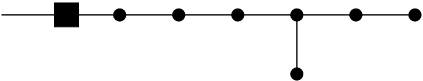}
\endrelabelbox} 
\medskip
\centerline{\it Figure 2: the $\lq$building block'.}
\medskip
\noindent Note that all circular dots here again represent tangent disc bundles. If $\tau \in \pi_{2k} BSO(2k)$ represents the tangent disc
bundle, then the square dot represents the bundle $3\tau.$
Plumbing according to the resulting graph yields the desired homotopy sphere as its boundary. 
For example, twice a generator ($2\tau$) is obtained from the following picture.
\bigskip
\centerline{
\relabelbox \epsfxsize 3 truein
\epsfbox{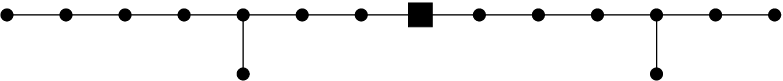}
\endrelabelbox} 
\medskip
\centerline{\it Figure 3: $2\tau \in \pi_{2k} BSO(2k)$.}
\medskip

Note that for each $\Sigma \in bP_{4k},$ [W1] establishes the existence (as opposed to the form) of corresponding simply-connected plumbing
diagrams with vertices representing $D^{2k}$-bundles over $S^{2k}.$ It is easily checked (for example following 
[Br;V.2]) that the graphs shown above have the required properties. However, to the best of the author's knowledge,
the explicit form of these graphs has not appeared in the literature before. 

As $bP_n$ is a cyclic group, it is clear that we can construct the homotopy sphere representing $q$-times the generator
corresponding to the $E_8$ graph by plumbing according to the graph consisting of $E_8$ and $(q-1 + p|bP_{4k}|)$ copies
of the $\lq$building block', for any integer $p \ge 0.$ In other words, for any fixed homotopy sphere $\Sigma \in bP_{4k}$, 
we can construct infinitely many different bounding manifolds via the plumbing construction. 
The technique of plumbing was originally introduced as a means of construcing manifolds with a given intersection form
on middle-dimension cohomology. In our situation, it is not difficult to see that for a homotopy sphere $\Sigma \in bP_{4k},$ the bounding
manifold discussed above has unimodular intersection form and signature $8(p|bP_{4k}|+q)$ (see
Propositions 1.4 and 1.5 in [W1]). Let us denote such a plumbed bounding manifold by $X_p$.

In [KS], Kreck and Stolz introduced the s-invariant of a spin manifold $M^{4k-1}$ which admits
a positive scalar curvature metric $g$ and has vanishing real Pontrjagin classes. This is a {\bb Q}-valued
invariant with the property that $|s|$ is constant on each component of the moduli space of positive scalar
curvature metrics, provided $H^1(M;\hbox{\bb Z}/2)=0$ (see [KS; Proposition 2.14]). The actual definition of 
the s-invariant is somewhat complicated and we refer
the reader to [KS] for details. However, for the manifolds under consideration in this paper
the s-invariant takes a simple form, as described in the following Lemma.
\proclaim Lemma 2.1. With $(M,g)$ as above, suppose that $W$ is a parallelisable manifold with $\partial W=M$,
equipped with a positive scalar curvature metric $\bar{g}$ which takes the form $dt^2+g$ in a neighbourhood 
of the boundary. Then the s-invariant
of $(M,g)$ is given by $$s(M,g)={1 \over {2^{2k+1}(2^{2k-1}-1)}} \sigma(W)$$ where $\sigma(W)$ is the signature of $W$.
\par
\smallskip
\noindent {\bf Proof:} The parallelisability of $W$ means in particular that $W$ is spin, so the following
formula [KS; Proposition 2.13] applies: $$s(M,g)=\hbox{ind}D^{+}(W,\bar{g})+t(W).$$
The positive scalar curvature metric on $W$ means that the term $\hbox{ind}D^{+}(W,\bar{g})$ vanishes
(see Remark 2.2(ii) in [KS]). The parallelisability of $W$ also means that the Pontrjagin classes
$p_1(W),...,p_n(W)$ vanish, and by equation (2.11) of [KS] this means that the term $t(W)$ is equal to
$1/(2^{2k+1}(2^{2k-1}-1)) \sigma(W).$ \hfill{\square}
\medskip
Note that the homotopy sphere $\Sigma \in bP_{4k}$ is a spin manifold with vanishing real Pontrjagin classes,
and that the bounding manifold $X_p$ is shown to be parallelisable in [W1].

\vfill\eject
\centerline{\bf \S3 Surgery}
\bigskip

Surgery is an operation which alters the topology of a manifold in a precise way, and was
introduced as a tool for use in problems concerning the classification of manifolds. A good introduction
to the subject is [R]. We describe the basic idea. 
Begin with a manifold $M^{n+m}$ and a smooth embedding $\iota_0:S^{n-1} \rightarrow M^{n+m}.$
Assume the normal bundle of $\iota_0(S^{n-1})$ is trivial. This means we can extend
$\iota_0$ to an embedding $\iota: S^{n-1} \times D^{m+1} \longrightarrow M.$
Performing surgery involves cutting something out - the interior of the image of the embedding $\iota$ -
and stitching something in, namely a copy of $D^n \times S^m$. 
Note that the boundary of $S^{n-1} \times D^{m+1}$ is $S^n \times S^{m-1}$, which agrees with the boundary of $D^n \times S^m$.
This gluing-in is performed in the obvious manner determined by $\iota$.
We say that the dimension of such a surgery is $n-1$, and the codimension is $m+1$.    
There is a subtlety in the above construction however: the extended embedding $\iota$ is not unique.
Viewing $S^{n-1} \times D^{m+1}$ as a bundle with base $S^{n-1}$ and fibre $D^{m+1}$, we can compose $\iota$ with a bundle
isomorphism $S^{n-1} \times D^{m+1} \rightarrow S^{n-1} \times D^{m+1}$ to obtain a new embedding $\phi.$ In general, performing
surgery using this new embedding will result in a manifold which is different from that obtained using $\iota$. Therefore
when performing surgery, it is imporant to specify the {\it trivialisation of the normal bundle}, that is,
the choice of embedding $\phi:S^{n-1} \times D^{m+1} \rightarrow M^{n+m}$ being used.   

As noted in \S2, any homotopy sphere of dimension $4k-1$ which bounds a parallelisable manifold arises as the boundary of a plumbed manifold,
where the disc bundles being plumbed are $D^{2k}$-bundles over $S^{2k}.$ 
(In fact, the same is true for homotopy spheres in dimensions $4k+1$ which bound a parallelisable manifold: see [LM; p162].)

It is an elementary observation that the effect on the boundary
of plumbing two such bundles $E_1$ and $E_2$ together is equivalent to performing surgery on a fibre $(2k-1)$-sphere in $\partial E_1$, using a trivialisation of the normal bundle determined by $E_2$ (see [W1; p.645]). Thus the homotopy spheres described in \S2 can be constructed, starting with an $S^{2k-1}$-bundle over $S^{2k}$, by performing successive surgeries. 
After performing a single surgery, the global $S^{2k-1}$-bundle structure will be destroyed. However, {\it
locally} this structure will be retained. Each of the successive surgeries required to construct the homotopy
sphere is a surgery on a $\lq$fibre' $S^{2k-1}$ arising from such a local bundle structure. In fact, this 
will always mean performing surgeries on the site of previous surgeries in an iterative manner. 

An approach to constructing Ricci positive metrics on homotopy spheres is to look for circumstances under which the surgery
operation can preserve Ricci positivity. The first major result in this direction was due to Sha and Yang [SY].
They assume that a tubular neighbourhood of the sphere on which surgery is to be performed
is isometric to the product of a round metric on the sphere and a round normal disk. Using this isometry as the
trivialisation of the normal bundle, they show that the manifold resulting from the surgery has positive
Ricci curvature provided the dimension of the surgery is at least one, the codimension at least three, and the ratio
of the radius of the sphere to the radius of the disk is sufficiently small. This result was developed further in [W2]
where it was shown that under the same metric assumptions, Ricci positivity can be preserved if different trivialisations
of the normal bundle are used. This generalisation comes at the expense of tighter dimensional requirements, but being
able to work with different trivialisations of the normal bundle is vital for applying the results to the construction of homotopy spheres.

The proof of these Ricci positive surgery results involves the construction of a Ricci positive metric on $D^n \times S^m$ which
glues smoothly with the metric on $M-\phi(S^{n-1} \times D^{m+1})$ when the surgery is completed. The metric on $D^n \times S^m$
is a {\it submersion metric}. The basic existence result for submersion metrics is due to Vilms [V, Theorem 9.59].

\noindent {\bf Construction 3.1.} {\it Let $G$ be a Lie group and $P \rightarrow B$ a principal $G$-bundle. Let 
$F$ be a manifold with a smooth (left) $G$-action and let $$E=P {\times}_G F
\longrightarrow B$$ be the associated fibre bundle over $B$ with fibre $F$. Suppose we 
have:
\item {(i)} a metric $\hat g$ on $F$ which is equivariant with respect to the
$G-$action;
\item {(ii)} a metric $\check g$ on $B$;
\item {(iii)} a principal connection $\nabla$ on P;
\item {(iv)} a smooth function $\mu: B \rightarrow \hbox{\bb R}^{+};$ \par
\noindent then there exists a unique submersion metric on $E$ with horizontal
distribution is associated to the connection, and given $b \in B$, $(F_b, g_b)$
is (non-canonically) isometric to $(F, \mu^2 \hat g)$.}
\smallskip 

\noindent In our situation, $E=D^n \times S^m$, which we will view as an $S^m$-bundle over $D^m,$ associated to the trivial $SO(m+1)$-principal
bundle over $D^n.$
The metric $\hat g$ on $S^m$ is just the round metric $ds^2_m.$ The metric $\check g$ is a warped product metric $dr^2+h^2(r)ds^2_{n-1}$
where $r$ is the radial parameter of the disc, and $h(r)$ is a smooth function. The fibre-scaling function $\mu:B \rightarrow \hbox{\bb R}^{+}$
will be a function $f(r).$ We will denote the resulting submersion metric on $E$ by $g(f,h,\nabla).$ The choice of connection $\nabla$
will depend on the particular normal bundle trivialisation being used in the surgery. In the case studied by Sha and Yang, a flat
connection sufficed, however, for homotopy sphere construction non-flat connections must be considered. Note that if $\nabla$ is flat,
$g(f,h,\nabla)$ is simply isometric to the warped product metric $dr^2+h^2(r)ds^2_{n-1}+ds^2_m.$ We will denote such a warped product
by $g(f,h).$ The principal connection $\nabla$ induces a horizontal distribution in $E.$ We will view $\nabla$ from now on as such
a horizontal distribution, to avoid having to invoke the principal bundle associated to $E.$

In the remainer of this section, we summarise the Ricci positive surgery results of [W2]. The statements we give below are slightly different to those presented in [W2], so as to emphasize the metric features which we will need in the sequel.

\noindent {\bf Theorem 3.2.} {\it Given $E=D^n \times S^m,$ a constant $\Delta \in (0,1)$ and a choice of connection $\nabla,$ there is a constant $\rho_0 \in (0,1),$
$\rho_0=\rho_0(n,m,\Delta,\nabla),$ and a constant $R' \in (0,1/4)$ independent of all parameters, such that for any choice of $\rho' \le \rho_0$ there exists a number $a=a(n,m,\Delta,\rho')>2$ depending smoothly on $\rho',$ with
$a \rightarrow \infty$ as $\rho' \rightarrow 0,$ 
and functions 
$f,h: [0,a] \rightarrow \bb R^{+}$ depending smoothly on $\rho'$ with the following properties:
\item{1)} $f(r) \equiv \rho'$ for all $r \in [0,R']$;
\item{2)} $f'(r) \ge 0$ and $f''(r) \ge 0$ for all $r$;
\item{3)} $f'(r)=\Delta$ in a small neighbourhood of $r=a$;
\item{4)} $h(r)=\sin r$ for all $r \in [0,R'];$ 
\item{5)} $h(r)$ is independent of all parameters (other than $r$) for $r \in [0,1/2]$;
\item{6)} $h(r)>0$ for $r>0$, $h'(r) \ge 0$ and $h''(r) \le 0$ for all $r$;
\item{7)} $h'(r) \equiv 0$ in a neighbourhood of $r=a$;
\item{8)} $\sup_{r \in [0,a]}|h''(r)|=\sup_{r \in [0,1/2]}|h''(r)|$;
\item{9)} $h(a)/f(a) \rightarrow 0$ as $\rho' \rightarrow 0$; \par
\noindent such that the submersion metric $g(f,h,\nabla)$ on $E$ has non-negative Ricci and positive scalar curvature, and strictly positive Ricci curvature
for $r$ small.} 
\smallskip

\noindent The functions $f$ and $h$ in this theorem have profiles as illustrated in Figure 4.
\smallskip

For the reader wishing to compare the above Theorem with the paper [W2], a few words of explanation are in order here.
The smooth dependence of $f,$ $h$ and $a$ on $\rho'$  might suggest that they are
all uniquely determined by $\rho'$ (given $n,$ $m$ and $\nabla$). However this is not the case, as certain choices have to
be made during the construction of these objects. On the other hand, it is clear from the construction details in [W2] that
the relevant choices can be made so that the variation with respect to $\rho'$ is smooth. 

Property (8) does not appear in [W2]
as it is not needed there. In fact, one can make choices in the construction of $h$ so that (8) fails, but all other properties
hold. Equally, one can make choices so that (8) holds. The argument is as follows: $h(r)$ is fixed for $r \in [0,1/2],$ and for $r>1/2$
the function is based on the solution to a certain ODE initial value problem, [W2; Definition 2.2 and Lemma 2.5]. The only difference between $h(r)$
for $r>1/2$ and this ODE solution occurs near $r=a.$ Here, the ODE solution has to be deformed in a concave down manner to satisfy property
(7) above (see [W2; Lemma 2.15]). This deformation takes place when $f'(r) \equiv \Delta,$ and as a result it is clear from the Ricci curvature
formulas [W2; equations (1.1)-(1.3)] that this deformation can take place over an interval of arbitrary length. By choosing this interval
to be suitably large (depending on the size of $\sup_{r \in [0,1/2]}|h''(r)|$), we can guarantee that (8) holds. This places a lower bound
on the size of $a,$ but has no other implications. 

The other properties of $f$ and $h$ listed in Theorem 3.2 are all evident from the construction
of these functions in [W2].   
\bigskip

\centerline{
\relabelbox \epsfxsize 4 truein
\epsfbox{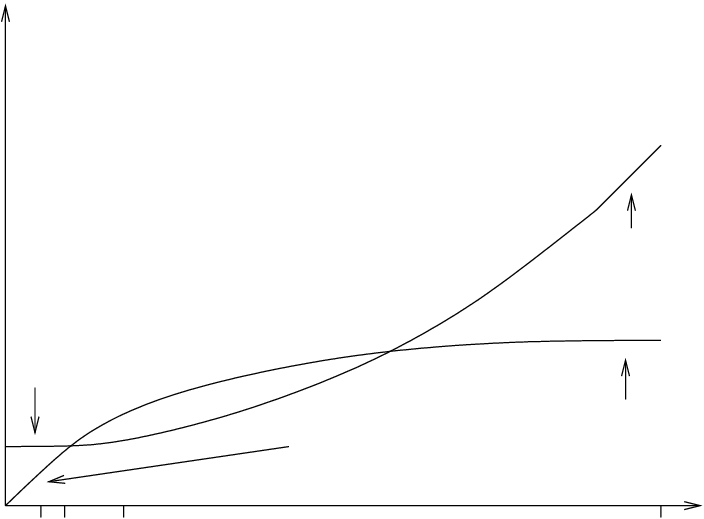} 
\extralabel <-14pt, 128pt> {$f(r)$}
\extralabel <0pt, 3pt> {$r$}
\extralabel <-14pt, 59pt> {$h(r)$}
\extralabel <-53pt, 88pt> {slope $\cos(R/N)$}
\extralabel <-238pt, 49pt> {slope 0}
\extralabel <-48pt, 32pt> {slope 0}
\extralabel <-140pt, 23pt> {$\sin r$}
\extralabel <-237pt, -10pt> {$R'$}
\extralabel <-225pt, -10pt> {${1 \over 4}$}
\extralabel <-205pt, -10pt> {${1 \over 2}$}
\extralabel <-20pt, -5pt> {$a$}
\extralabel <-255pt, 23pt> {$\rho'$}
\endrelabelbox}
\bigskip
\centerline{\it Figure 4}
\bigskip

The above existence result can easily be adapted for Ricci positive surgery. 
Suppose we have the following:
\smallskip
\item{(i)} A Riemannian manifold $(M^{n+m},g)$ with positive
Ricci curvature.
\smallskip
\item{(ii)} An isometric embedding $\iota :S^{n-1}(\rho) \times D^{m+1}_R(N) \rightarrow
M$ where $D^{m+1}_R(N)$ denotes a geodesic ball of radius $R$ in the $(m+1)$-sphere with the
round metric of radius $N$. $S^{n-1}(\rho)$ is the $(n-1)$-sphere with the round
metric of radius $\rho$.
\smallskip
\item{(iii)} A smooth map $T:S^{n-1} \rightarrow SO(m+1)$. This
induces a map $$\tilde T:S^{n-1} \times D^{m+1} \rightarrow S^{n-1}
\times D^{m+1}$$
$$(x,r,y) \mapsto (x,r,T(x)y)$$
where we are using polar coordinates $(r,y)$ in $D^{m+1},$ with $r \in [0,R]$ and 
$y \in S^m$. In turn, this induces a map $\phi:S^{n-1} \times D^{m+1} \rightarrow M$
where $\phi:=\tilde T \circ \iota.$
\smallskip
\noindent Define the manifold $\hat M$ to be
$$ \hat M = (M \setminus \phi (S^{n-1}(\rho) \times D^{m+1}_R(N)) \amalg
(D^n \times S^m)/ \sim$$
where $\sim$ indicates
the identification of boundaries using $\phi$ in the obvious manner.
\smallskip
We now come to the main surgery result in [W2]. To state this, we will denote by $a_0$ the value of $a$ in Theorem 3.2 corresponding
to $\rho'=\rho_0.$ Similarly, the corresponding functions $f$ and $h$ will be denoted $f_0,$ $h_0.$ 
\medskip

\noindent{\bf Theorem 3.3.} {\it Suppose we have $(M^{n+m},g)$, $\rho,$ $R,$ $N$ and $\phi$ as specified in $(i)-(iii)$ above. Suppose further
that $n \ge m+1 \ge 3$. Set $\Delta=\cos(R/N),$ and choose $\delta \in (0,R')$ where $R'$ is the universal
constant appearing in Theorem 3.2. Then there is a connection $\nabla$ on $E=D^n \times S^m$ which is flat outside the
region $r \in (\delta/3,2\delta/3),$ so that provided
$$\rho<N\sin(R/N) {{h_0(a_0)} \over {f_0(a_0)}} \eqno{(\dag)}$$
is satisfied, for any choice of $\rho' \le \rho_0,$  
the Ricci non-negative submersion metric $g(f,h,\nabla)$ of Theorem 3.2 can be made to
fit smoothly into
$$M \setminus \phi (D^{m+1}_R(N) \times S^{n-1}(\rho))$$ 
by
\smallskip
\item{(a)} globally rescaling $g(f,h,\nabla)$ by a factor of $N\sin(R/N)/f(a),$ and then 
\smallskip
\item{(b)} deforming the rescaled $f(r)$ to join smoothly with the function $N\sin((R+r-a)/N)$ at $r=a.$
\smallskip
\noindent Moreover, this deformation can be performed in such a way that the resulting metric on $\hat{M}$ has non-negative Ricci curvature.}
\bigskip

\noindent {\bf Remark 3.4.} The metric on $\hat{M}$ in Theorem 3.3 has strictly positive Ricci curvature on $M-\phi(D^{m+1} \times S^{n-1}) \subset
\hat{M},$ as well as in the central region of the $\lq$glue-in' piece $E.$ The fact that the metric elsewhere might
only be Ricci non-negative is, however, not a concern. By a result of Ehrlich ([E], but see also [A]), if a Ricci non-negative
manifold has strictly positive Ricci curvature at some point, the positive curvature can be spread out from that point via a smooth
metric deformation to cover the whole manifold. We are therefore justified in viewing Theorem 3.3 as a Ricci positive surgery
theorem: the Ricci non-negative Riemannian manifolds we create using the Theorem can, as a final step, be smoothly defomed to exhibit
Ricci positivity. (In fact, it is not difficult to see how to obtain such a Ricci positive metric explicitly: simply adjust the
functions $f(r)$ and $h(r)$ by an arbitrarily small amount in neighbourhoods with some zero Ricci curvatures, so as to create strict
concavity at all values of $r.$)
\smallskip 

Suppose we wish to perform the kind of metric surgery described by Theorem 3.3 on the fibre-sphere of a sphere bundle over a Ricci positive
base manifold. It is well-known (see [Be; 9.70]) that the submersion metric on such a bundle constructed by using the given Ricci positive
base metric, scaled round metrics $\rho ds^2$ for some constant $\rho>0$ on the fibres, and any choice of connection,
will have positive Ricci curvature provided $\rho$ is chosen sufficiently small. If we assume that the base metric in a neighbourhood
of the surgery sphere is isometric to $D^{m+1}_R(N)$ for some $R$ and $N$ (and we can always arrange 
for this to be the case via a defomation of the base metric), then Theorem 3.3 applies. For any choice of trivialisation $\phi$ we obtain a Ricci positive metric on the manifold $\hat M$ resulting from the surgery, as it is clear that we can choose $\rho$ as small as we like, and so can guarantee
that condition $(\dag)$ in Theorem 3.3 is satisfied. 

Properties (1) and (4) in Theorem 3.2, together with the fact that we are free in Theorem 3.2 to choose $\rho'$ to be arbitrarily small,
means that the metric on $\hat M$ resulting from Theorem 3.3 is of the correct form to attempt further
surgeries using Theorem 3.3 on the $S^m$-fibres of the glue-in piece $E=D^n \times S^m.$ As noted above, an iterative sequence of such
surgeries are required to construct the homotopy spheres under consideration in this paper. In particular, repeated application
of Theorem 3.3 shows that all homotopy spheres which bound parallelisable manifolds admit Ricci positive metrics. 

\proclaim Convention 3.5. In order to use Threorem 3.3 to perform a metric surgery on the site of a previous surgery (as discussed in the above
paragraph), we will always set $R=R'/4,$ $N=1$ and $\rho=\rho',$ where $\rho'$ is the quantity arising from the previous surgery.
\par

\noindent The significance of this convention will be made clear \S4.

\proclaim Notation 3.6. Let $\Sigma^{4k-1}$ be a homotopy sphere which bounds a parallelisable manifold. 
As described in \S2, view $\Sigma$ as the boundary of the plumbed manifold $X_p.$ Denote by $g_p$ a 
Ricci positive metric on $\Sigma$ obtained from a Ricci positive submersion metric on the tangent sphere bundle of $S^{2k}$ by repeated
Ricci positive surgeries of the form described in Theorem 3.3, performed according to the plumbing diagram for $X_p.$ \par

\bigskip\bigskip
\centerline{\bf \S4 A metric deformation}
\bigskip
As described in \S1, given a homotopy sphere $\Sigma$ which bounds a parallelisable manifold equipped with a Ricci positive $\lq$plumbed' metric $g_p$ (see 3.6), to prove Theorem A we need to show that $g_p$ can be extended to a positive scalar curvature metric on the plumbed bounding manifold $X_p$ (see \S2) in such a way that it is a product near the boundary. Performing any kind of positive scalar
curvature extension of $g_p$ across $X_p$ turns out to be quite difficult. The situation can be made more tractable, however, by first
{\it deforming} $g_p$ smoothly through positive scalar curvature metrics on $\Sigma$ to a new metric $g_{p,\infty}.$ This metric does not have
positive Ricci curvature, but it belongs to the same component of the space of positive scalar curvature metrics as $g_p,$ and this
is sufficient for our purposes. Because of the form of $g_{p,\infty},$ it is straightforward to extend to a positive scalar curvature metric on $X_p.$
This extension is carried out in \S5. Unfortunately, this extended metric is not a product in a neighbourhood of the boundary. However, because
of the form of $g_{p,\infty},$ this is easy to correct via a metric deformation within positive scalar curvature. This is carried out in \S6. Thus the Kreck-Stolz invariant for $(\Sigma,g_{p,\infty})$ can be computed using Lemma 2.1 in terms of the
signature of $X_p,$ leading to a proof of Theorem A. 

In the current section, we show how to deform $g_p$ to $g_{p,\infty}.$

The metric $g_p$ naturally splits into a number of $\lq$pieces' arising from the various surgeries used in its construction.
Apart from the initial tangent sphere bundle submersion metric, these pieces are all formed 
from the metrics $g(f,h,\nabla)$ in Theorem 3.2 by global rescale and boundary adjustment (see Theorem 3.3). 
Our first goal is to describe a deformation of one of these pieces, and
our starting point is to consider the simplest case of a metric $g(f,h,\nabla):$ the case where the connection $\nabla$ is flat.
In this case, $g(f,h,\nabla)$ is isometric to the warped product $g(f,h)=dr^2+h^2(r)ds^2_n+f^2(r)ds^2_m.$ 

\proclaim Lemma 4.1. Consider the metric $g(f,h)=dr^2+h^2(r)ds^2_n+f^2(r)ds^2_m$ on $D^n \times S^m.$ There exists a constant $c_0=c_0(n)$  
such that for all $0<c<c_0$ the metric $g(f,h)$ can be deformed smoothly through positive scalar curvature metrics to a metric 
$g(f,h_{\infty}),$ where $h_{\infty}(r)$ is a concave down function satisfying
$h_{\infty}(r)=\sin r$ for $r \in [0,\delta]$ where $\delta=\delta(c)<R'$, and $h_{\infty}(r) \equiv c$ for $r>R'$. (See Figure 5
for an illustration of this function, and Theorem 3.2
for the definition of $R'$). Moreover, this deformation can be performed keeping the metric constant for $r \in [0,\delta].$
\par

\noindent{\bf Proof.} The deformation will be parameterised using a $\lq$time' parameter $\tau \in [0,1].$ 
Set $h_{\tau}:=(1-\tau)h+\tau h_{\infty}.$ Notice that $h_{\tau}$ is a concave down function for each $\tau.$
We need to study the conditions under which the family of metrics $g(f,h_{\tau})$ have positive scalar curvature for each $\tau.$

An elementary calculation shows that the scalar curvature of this metric is given by
$$-2(n-1){{h_{\tau}''} \over {h_{\tau}}}-2m{{f''} \over f}+(n-1)(n-2){{1-(h_{\tau}')^2} \over {h_{\tau}^2}}+m(m-1){{1-(f')^2} \over {f^2}}-2(n-1)m{{h_{\tau}'f'} \over {h_{\tau}f}}$$ where all the derivatives are taken with respect to $r.$

First note that there is nothing to prove for $r \in [0,R']$ as $f(r)$ is constant at these values of $r,$ and by concavity, the sum of the terms involving only $h_{\tau}$ and its derivatives is strictly positive for all $r$ and for each $\tau \in [0,1].$ From now on we therefore assume that $r \ge R'.$ 

Next, we claim that the term $h_{\tau}'f'/h_{\tau}f$ is decreasing with respect to $\tau$ for each $r.$
As $f$ is independent ot $\tau,$ it suffices to consider the behaviour of $$\eqalign{
{{h_{\tau}'} \over {h_{\tau}}}&={{(1-\tau)h'+\tau h'_{\infty}}
\over {(1-\tau)h+\tau h_{\infty}}} \cr &={{h'} \over {h+{{\tau} \over {1-\tau}}c}} \cr
}$$ for these values of $r$.
Now $\tau/(1-\tau)$ is strictly increasing for $\tau \in [0,1),$ and this establishes the claim.

\bigskip\medskip
\centerline{
\relabelbox
\epsfbox{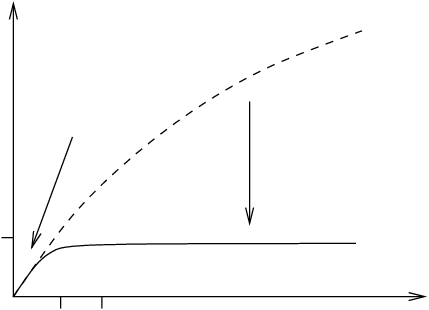}
\extralabel <-30pt, 30pt> {$h_\infty(r)$}
\extralabel <-0pt, 5pt> {$r$}
\extralabel <-30pt, 135pt> {$h(r)$}
\extralabel <-183pt, -10pt> {$\delta$}
\extralabel <-165pt, -10pt> {$R'$}
\extralabel <-175pt, 87pt> {$\sin(r)$}
\extralabel <-220pt, 34pt> {$c$}
\extralabel <-85pt, 71pt> {deform}
\endrelabelbox}
\bigskip
\centerline{\it Figure 5}
\bigskip

To complete the proof of Lemma 4.1, it therefore suffices to check that $-2h''_{\tau}/h_{\tau}+(n-2)(1-(h_{\tau}')^2)/h_{\tau}^2$ is non-decreasing with $\tau$ for all $r \in [R',a].$ 
To this end, a straightforward calculation shows that the derivative of this expression with respect to $\tau$ is equal to 
$${2 \over {[(1-\tau)h+\tau c]^3}} \Bigl[ ch''[(1-\tau)h+\tau c]+(n-2)(h-c)+(n-2)c(1-\tau)(h')^2 \Bigr]. $$
As $2[(1-\tau)h+\tau c]^{-3}>0$, we only need check the positivity of the term $$ch''[(1-\tau)h+\tau c]+(n-2)(h-c)+(n-2)c(1-\tau)(h')^2.$$ 
By Theorem 3.2 (6) we see that the last term in this expression is always non-negative, the middle term is strictly positive
if $c<h(R')=\sin(R')$ as $h$ is increasing, but the first term is non-positive as $h'' \le 0.$ By 3.2 (8)
we have that the absolute value of the first term is bounded above by $c(h(r)+\tau c)\sup_{r \in [0,1/2]}|h''(r)|.$
For simplicity, let us set $\lambda:=\sup_{r \in [0,1/2]}|h''(r)|.$ Thus it suffices to show that $$c\lambda[(1-\tau)h(r)+\tau c]<(n-2)(h(r)-c)$$
for all $r \in [R',a].$ Rearranging this we obtain $$(n-2)c+\tau c^2 \lambda < h(r)[n-2-c\lambda(1-\tau)].$$ As $h$ is increasing and $\tau \in [0,1]$
it suffices to establish $$(n-2)cf+c^2\lambda < \sin(R')[n-2-c\lambda],$$ which is clearly true for all $c$ sufficiently small. Moreover, the
upper bound $c_0$ for the values of $c$ which satisfy this inequality clearly depends only on $n,$ $\lambda$ and $R'.$ But $\lambda$ and $R'$
are constants independent of any other quantities, so we are therefore justified in writing $c_0=c_0(n).$

Having fixed a value for $c$, we finally set $\delta$ to be any number consistent with the requirements for
$h_{\infty}$ specified in the statement of the Lemma. 
\hfill{\square}
\medskip

So far we have only considered deformation of the warped product metric $g(f,h).$ It is now time to bring connections back into the picture.
\proclaim Proposition 4.2. Suppose we have $(M^{n+m},g)$, $\rho,$ $R,$ $N,$ $\phi$ and $\Delta$ as in Theorem 3.3.
Fix a value of $c<c_0(n),$ where $c_0(n)$ is the constant appearing in Lemma 4.1. Set $\delta =\delta(c)\in (0,R')$ as in Lemma 4.1.
Then there is a connection $\nabla$ on $E=D^n \times S^m$ which is flat outside the
region $r \in (\delta/3,2\delta/3),$ so that provided condition $(\dag)$ in Theorem 3.3 is satisfied, the Ricci non-negative submersion 
metric $g(f,h,\nabla)$ admits a smooth deformation through positive scalar curvature metrics to $g(f,h_{\infty},\nabla),$ with $h_{\infty}$ as in Lemma 4.1. Moreover, these metrics (after rescale and boundary adjustment) give a smooth family of positive scalar curvature metrics on $\hat{M},$ the manifold resulting from
the surgery. \par
\medskip
\noindent {\bf Proof.} It sufficies to note that since $(\delta/3,2\delta/3) \subset [0,R'],$ in the region where the connection is non-flat
the path of metrics is fixed independently
of $\tau.$ On the other hand, for $r>R',$ the metric $g(f,h,\nabla)$ is isometric to the warped product $g(f,h),$
for which the result was established in Lemma 4.1. The final statement is obvious.
\hfill{\square}
\medskip

We immediately deduce:
\proclaim Corollary 4.3. For any collection of Ricci positive surgeries of the form described by Theorem 3.3, uniform values for 
$c$ and $\delta$ can be chosen in advance.
\par

\proclaim Convention 4.4. For any sequence of metric surgeries we wish to perform, we will assume, in accordance with Corollary 4.3, that
we use the same fixed values for $c$ and $\delta$ for all surgeries. Moreover, we will also assume (as we did in Proposition 4.2) that
the value of $\delta$ in Theorem 3.3 (used to bound any region of non-flatness for the connection) is precisely $\delta(c).$ \par

The significance of Convention 4.4 is that by constraining the non-flat part of the connection to lie in the region $r \in (\delta(c)/3,
2\delta(c)/3),$ the metric deformations being described in this section take place where the metric is a warped product, as do subsequent surgeries, which by convention Convention 3.5 take place at $r \le R'/4 \, (<\delta(c)/3).$ 

\smallskip

Having discussed deforming an individual surgery $\lq$piece' through positive scalar curvature metrics, we now want to show how this
procedure can be iterated to deform the metric $g_p$ to the desired metric $g_{p,\infty}.$

For a given homotopy sphere $\Sigma \in bP_{4k},$ consider the corresponding plumbing diagram, for example one of those illustrated in
Figure 1 or Figure 3.

To interpret this diagram as a template for surgery, let us read it from right to left, thus the right-most vertex will represent the
initial sphere bundle, and each of the vertices to its left represents a successive surgery. All such graphs bifurcate in a number of places,
as illustrated in Figure 6.
\bigskip
\centerline{
\relabelbox \epsfxsize 0.5 truein
\epsfbox{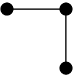}
\endrelabelbox} 
\medskip
\centerline{\it Figure 6}
\medskip
In the following we will assume that both the original metrics and the subsequent deformation operations we perform
on the objects represented by the vertex to the left
and the vertex below the bifurcation point are identical. As a result, it suffices to consider {\it linear} plumbing graphs.

Suppose we have a linear graph with $i_0+1$ vertices, and wish to use this as a surgery template.
Assume the right-hand vertex is a sphere bundle over a sphere, equipped with a Ricci positive submersion metric having round fibres of radius $\rho$ and a round base
sphere of radius $N.$ Without loss of generality, we can assume that the connection on this bundle is chosen to be flat in a neighbouhood
of the surgery, and that the base disc removed has radius $R.$ 
Assume further that the dimension of each surgery is $\ge$ codimension $\ge 3,$ so Theorem 3.3 applies.

Let us introduce some notation to help clarify later constructions. Let $\rho$ (as before) be the radius of the fibre sphere for the right-hand vertex,
and let $\rho'_i$ be the value of $\rho'$ from Theorem 3.3 for the $i^{th}$ surgery, that is, for the $(i+1)^{st}$ vertex from the right.
Label all other quantites associated with this surgery by a subscript $i.$ Let $M$ be the manifold resulting from the
surgery process, and let $g$ be the Ricci non-negative metric on $M$ obtained from successively applying Theorem 3.3. 

\proclaim Lemma 4.5. For each $t \in [0,1]$ there is a Ricci non-negative metric $g_t$ on the manifold $M$ (as defined above) constructed via metric surgery in the same way as $g,$ such that the $\rho'$ value for the left-most vertex is $t\rho'_{i_0},$ and the family $g_t$ depends smoothly on $t$.
\par
\medskip
\noindent {\bf Proof.} We proceed by induction on the number of vertices. Consider first the case of a single surgery. The corresponding graph consists
of two vertices, as shown in Figure 7.
\bigskip
\centerline{
\relabelbox \epsfxsize 0.5 truein
\epsfbox{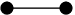}
\endrelabelbox} 
\medskip
\centerline{\it Figure 7}
\medskip
As discussed above, the right-hand vertex represents a sphere bundle with Ricci positive submersion metric.
The left-hand vertex represents a certain surgery on the right hand vertex, where the metric on the glue-in piece is the metric from Theorem 3.3. As pointed out in Theorem 3.3, this metric will differ from the metric $g_1(f_1,h_1,\nabla_1)$ by a global
rescale and a concave down adjustment to the (rescaled) function $f_1$ near the boundary. 

Let us consider the metric $g_1(f_1,h_1,\nabla_1).$ Replacing $\rho'_1$ by $t\rho'_1$ for decreasing values of $t<1$ results in increasing values for $a_1,$ with the dependence on $t$ being smooth (see 3.2). In turn, this forces a smoothly decreasing global rescale factor, $\lambda_1(t)$ say, in order to complete the surgery.
The rescale is necessary in order for $f_1$ to join smoothly with the fuction $N_1\sin((R_1-a_1+r)/N_1)$ at $r=a_1.$ 
However, the smaller rescale factor also has the effect of lowering the (rescaled) value of $h_1(a_1)$ below $\rho.$ This is not a problem as we are free to smoothly reduce $\rho$ accordingly, that is, to shrink the radii of the sphere bundle fibres. This completes the
initial step of the induction argument.

Now suppose the result is true for a linear graph of $i_0-1$ vertices. Consider the $i_0^{th}$ vertex and the metric 
$g_{i_0}(f_{i_0},h_{i_0},\nabla_{i_0}).$
Replacing $\rho'_{i_0}$ by $t\rho'_{i_0}$ for decreasing values of $t<1$ forces a smaller global rescale factor $\lambda_{i_0}(t)$ for $g_{i_0}(f_{i_0},h_{i_0},\nabla_{i_0}).$ Just as
for the initial case, this necessitates correspondingly smaller values for $\rho'_{i_0-1}.$ But by our inductive assumption, we can re-choose
$\rho'_{i_0-1}$ in this way while inducing smoothly varying Ricci non-negative metrics on the manifold corresponding to the first
$i_0-1$ vertices of the graph. This completes he inductive step.
\hfill\square

\proclaim Lemma 4.6. Consider again $(M,g)$ as in Lemma 4.5. For each $\tau \in [0,1],$ there is a positive scalar curvature metric $g_{\tau}$ on $M$ constructed via metric surgery in the same way as $g,$ such that the metric on the leftmost vertex corresponds to 
$g_{i_0}(f_{i_0},(h_{i_0})_{\tau},\nabla_{i_0})$ and the family $g_{\tau}$ depends smoothly on $\tau$.  
\par
\smallskip
\noindent {\bf Proof.} Consider the metric $g_{i_0}(f_{i_0},h_{i_0},\nabla_{i_0}).$ By re-choosing the connection $\nabla_{i_0}$ if necessary so as to be flat outside the region $r \in (\delta/3,2\delta/3),$ where $\delta$ is determined by any fixed choice of $c \in (0,c_0)$ as in Lemma 4.1, we see
from Proposition 4.2 that the family of metrics $g_{i_0}(f_{i_0},(h_{i_0})_{\tau},\nabla_{i_0})$ has positive scalar curvature. Each of these
metrics will glue smoothly with the metric to its $\lq$right' in the graph provided $\rho'_{i_0-1}$ is re-chosen to agree with the final
(rescaled) value of $(h_{i_0})_{\tau}.$ See Figure 8 below. The result now follows from Lemma 4.5. 
\hfill\square 
\smallskip

\font\eightrm=cmr8
\centerline{
\relabelbox \epsfxsize 4.5 truein
\epsfbox{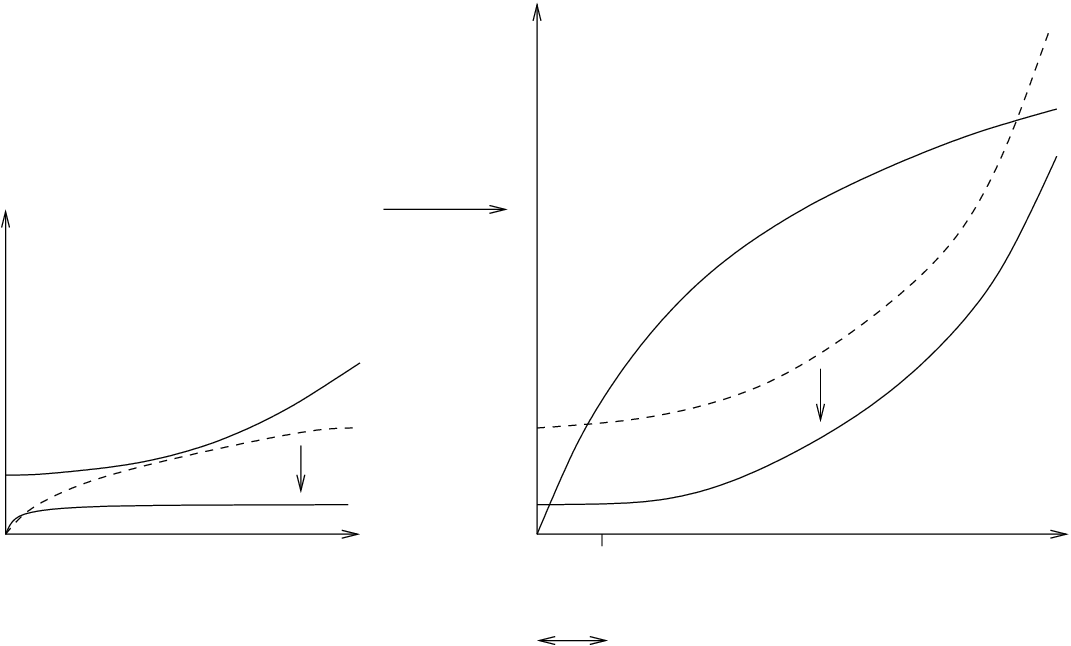} 
\extralabel <-92pt, 58pt> {\sevenrm adjust}
\extralabel <-180pt, 24pt> {$r$}
\extralabel <0pt, 26pt> {$r$}
\extralabel <-95pt, 108pt> {$h_1$}
\extralabel <-20pt, 84pt> {$f_1$}
\extralabel <-31pt, 159pt> {\sevenrm [original $\scriptstyle f_1$]}
\extralabel <-126pt, 14pt> {${{R'} \over 4}$}
\extralabel <-143pt, -5pt> {\eightrm remove}
\extralabel <-150pt, -11pt> {\eightrm for surgery}
\extralabel <-150pt, 35pt> {$\rho'_1$}
\extralabel <-153pt, 55pt> {$\scriptstyle [\rho'_1]$}
\extralabel <-195pt, 44pt> {\sevenrm deform}
\extralabel <-183pt, 34pt> {$\scriptstyle h_{2,\infty}$}
\extralabel <-180pt, 54pt> {$\scriptstyle h_2$}
\extralabel <-181pt, 70pt> {$\scriptstyle f_2$}
\extralabel <-171pt, 123pt> {\eightrm push}
\extralabel <-177pt, 115pt> {\eightrm together}
\extralabel <-290pt, 42pt> {$\scriptstyle \rho'_2$}
\endrelabelbox}
\bigskip\medskip\vskip 1pt
\centerline{\it Figure 8: metric deformation in the second surgery $\lq$piece'.}
\bigskip
\proclaim Observation 4.7. The deformation of Lemma 4.5 can be applied to any of the surgery $\lq$pieces' which make up $M,$ as it has
no effect on that part of the manifold represented by vertices to the left.
\par

\proclaim Corollary 4.8. The deformations can be performed successively from the left of the diagram to create an overall deformation from $g_p$ to a new
positive scalar curvature metric $g_{p,\infty}.$
\par

\bigskip\bigskip
\centerline{\bf \S5 Handles}
\bigskip

Recall from \S3 that given an embedding $\phi:S^{n-1} \times D^{m+1} \rightarrow M,$ we can perform surgery on $M,$ which results in a new
manifold $\hat{M}$ given by 
$$ \hat M = (M \setminus \phi (S^{n-1}(\rho) \times D^{m+1}_R(N)) \amalg_{\phi}
(D^n \times S^m).$$
Now suppose that $M=\partial {\bar{M}}.$ We can perform a related construction on $\bar{M}$ using $\phi.$ Note that $$\partial (D^n \times D^{m+1})
=(S^{n-1} \times D^{m+1}) \cup (D^n \times S^m).$$ Using the map $\phi$ on the first $\lq$half' of this boundary, we can glue $D^n \times D^{m+1}$
onto $\bar{M},$ to yield a new manifold with boundary $(D^n \times D^{m+1}) \cup_{\phi} \bar{M}.$ Notice that on the boundary, this construction
is just the original surgery. Here, $D^n \times D^{m+1}$ is called a {\it handle} (or handlebody), the process just described is called {\it handle
addition.} 

Notice that the plumbing of two disc-bundles over spheres can be viewed as attaching a
handle to either of the bundles. This means, for example, that the plumbed manifolds described in \S2 having homotopy spheres as boundary can be viewed as constructed from a disc bundle over a sphere by successively adding handles.

Our main aim in this section is to prove a result about extending surgery metrics across the corresponding handles (Proposition 5.2). 
However, we first need a lemma concerning metrics on disc bundles.
We will say that a smooth Riemannian manifold with boundary 
has positive scalar (respectively positive Ricci) curvature when the metrics on the interior and the boundary both have positive
scalar (respectively Ricci) curvature.  
\medskip
\proclaim Lemma 5.1. Suppose that $E$ is the total space of a bundle with base $M$, a manifold with a fixed Ricci positive metric, fibre $S^{n-1}$, and structural group
$SO(n)$. Let $\theta(s)$, $s \ge 0$, be any function satisfying: $\theta$ is odd at $s=0$, $\theta'(0)=1$,
$\theta''(s)<0$ for all $s>0$, $\theta'(s)>0$ for all $s$, and $\lim_{s \rightarrow 0^{+}} \theta''(s)/\theta(s)<0.$
For a given choice of principal connection $\nabla$ on the principal bundle associated to $E$, there exists
$R_0=R_0(\nabla,\theta)$ such that for all $\rho \le R_0$, the submersion metric on $E$ with fibres isometric to $\theta^2(\rho) ds^2_{n-1}$, and the submersion metric
on the corresponding $D^n$-bundle $\bar{E}$ with fibres isometric to $ds^2+\theta^2(s)ds^2_{n-1}$ where $s \in [0,\rho],$ both have positive Ricci curvature. \par
\medskip
\noindent {\bf Proof:} We view the metric on $\bar{E}$ as $ds^2+g_s$, where $g_s$ is the induced submersion metric on $E$ with
fibre radius $\theta(s)$. The following Ricci curvature formulas for $ds^2+g_s$ are easily computed (for example
by adapting the formulas in [W3; Proposition 4.2]).
$$\eqalign{ \hbox{Ric}&(\partial/\partial s,\partial/\partial s)=-(n-1){{\theta''}\over {\theta}} \,; \cr 
\hbox{Ric}&(X_i,X_i)=\hbox{Ric}_M (X_i,X_i)-2\theta^2\langle A_{X_i},A_{X_i} \rangle \,; \cr 
\hbox{Ric}&(U_j,U_j)={{(n-2)} \over {\theta^2}} \, (1-\theta'^2)-{{\theta''} \over {\theta}}
+\theta^2 \langle AV_j,AV_j \rangle \, ; \cr
\hbox{Ric}&(X_i,U_j)=-\theta\langle (\delta\check{A})X_i,V_j \rangle. \cr }$$
Here, $\{ X_i \}$ are (the lifts to $E$ of) orthonormal vector fields tangent to the base $M$; $\{ U_j \}$ are 
orthonormal vector fields tangent to the distance spheres $S^{n-1} \subset D^n$, and the $\{ V_j \}$ are the corresponding vector
fields which would be orthonormal with respect to the metric $ds^2_{n-1}$; the $A$-tensor expressions appearing here are
those for the unit sphere bundle diffeomorphic to $E$ (see [Be; \S9C] for definitions), 
and so are independent of $s$. 

The concavity of $\theta$ ensures that $\hbox{Ric}(\partial/\partial s,\partial/\partial s)>0.$
The Ricci positivity of $M$ ensures $\hbox{Ric}(X_i,X_i)>0$ provided $s$ (and hence $\theta(s)$) is suitably small.
Similarly, the positivity of $\hbox{Ric}(U_j,U_j)$ is clear. 

It remains to check the influence of the $\lq$mixed' term $\hbox{Ric}(X_i,U_j)$. 
(Note that the other mixed terms $\hbox{Ric}(\partial/\partial s,X_i)$ and $\hbox{Ric}(\partial/\partial s,U_j)$ vanish.) To ensure that this term does
not upset Ricci positivity we need to check that $$\hbox{Ric}(X_i,X_i) \hbox{Ric}(U_j,U_j)>(\hbox{Ric}(X_i,U_j))^2.$$
Using the above expressions and considering both sides of this inequality as $s \rightarrow 0$, we arrive at the limiting
inequality $$\hbox{Ric}_M(X_i,X_i)\Bigl(\, \lim_{s \to 0^{+}}\, {{(n-2)} \over {\theta^2(s)}} \, (1-\theta'^2(s))-{{\theta''(s)} \over {\theta(s)}}\Bigr)>0,$$ which is true. Thus the bundle metric
on $\bar{E}$ is Ricci positive for $\rho$ suitably small. It is also clear (for example by [Be; \S9G]) that the
metrics $g_s$ on $E$ have positive Ricci curvature for all suitably small $s$. \hfill{\square}
\medskip 

\proclaim Proposition 5.2. Suppose we have $(M^{n+m},g)$, $\rho,$ $R,$ $N,$ $\phi$ and $\Delta$ as in Theorem 3.3. Suppose further that
$M=\partial \bar{M},$ and that the metric $g$ extends to a positive scalar curvature metric $\bar{g}$ on $\bar{M}.$ 
Using the outward normal parameter to the boundary, extend $\phi$ to an embedding $\bar{\phi}:[-\epsilon, 0] \times S^{n-1} \times D^{m+1}
\rightarrow \bar{M}.$ Assume that the metric on $\bar{M}$
pulled back to this product takes the form $$ds^2+\zeta^2(s)ds^2_{n-1}+dt^2+N^2\sin^2(t/N)ds^2_m,$$
where $s$ parametrises $[-\epsilon, 0]$ and $t \in [0,R]$ is the radial parameter of $D^{m+1}$. 
In particular, this means that $\zeta(0)=\rho.$ Suppose further that
$\zeta'(s)>0$ and $\zeta''(s)<0$ for all $s \in [-\epsilon, 0]$. 
Then there is a connection $\nabla$ on $E=D^n \times S^m$ which is flat outside the
region $r \in (\delta/3,2\delta/3),$ so that assuming condition $(\dag)$ in Theorem 3.3 holds, the positive scalar curvature metric $g(f,h_{\infty},
\nabla)$ on $E$ from Proposition 4.2, after globally rescaling and a suitable adjustment to $f$ near the boundary, extends to 
a positive scalar curvature metric on $\bar{E}=D^n \times D^{m+1},$ which glues smoothly via $\phi$ to give a positive scalar curvature
metric on $\bar{M} \cup_{\phi} \bar{E}.$
\par
\medskip
\noindent {\bf Proof.} We will first consider the situation where the connection $\nabla$ can be chosen to be flat. In this case, the
metric $g(f,h_{\infty},\nabla)=g(f,h_{\infty})$ on $E$ is a warped product metric. Now $g(f,h_{\infty})$ has to be globally rescaled
and then $f$ adjusted near the boundary of $E$ in order to use it for positive scalar curvature surgery on $M$ (see the comments in Theorem 3.3).

Note that if the global rescale factor is $\lambda^2,$ then the rescaled metric is $$\lambda^2 dr^2+ \lambda^2 h^2_{\infty}(r/\lambda)
ds^2_{n-1}+\lambda^2 f^2(r/\lambda) ds^2_m$$ for $r \in [0,a].$
By a parameter change $\tilde{r}:=\lambda r,$ this is equivalent to $$d\tilde{r}^2+\lambda^2 h_{\infty}^2(\tilde{r}/\lambda)ds^2_{n-1}+\lambda^2f^2(
\tilde{r}/\lambda)ds^2_m$$
for $\tilde{r} \in [0,\lambda a].$ That is, by introducing a rescaled parameter, we can write the rescaled metric in the form $g(\tilde{f},{\tilde{h}}_{\infty})$
for $\tilde{r} \in [0,\tilde{a}],$ with $\tilde{f}=\lambda f,$ ${\tilde{h}}_{\infty}=\lambda h_{\infty}$ and $\tilde{a}=\lambda a.$ 
In the interests of notational simplicity, we will assume from now on,
by ignoring the tilde, that the parameter $r$, the number $a$ and the functions $f(r)$ and $h_{\infty}(r)$ all refer to
the relevant quantities for the {\it rescaled} metric.

The adjustment of $f$ to meet the boundary conditions required for smooth gluing (that is, the bending of $f$ from having constant
slope near the boundary to joining smoothly with $N\sin((R-a+r)/N)$ at $r=a$) will need our careful attention. In particular, we claim that
this adjustment can be made within positive scalar (and in fact locally, non-negative Ricci) curvature in such a way that the inequality $$f'(r) \le \sqrt{
1-(f(r)/N)^2} \eqno{(\ast)}$$ holds for all $r.$ This inequality will be important later on in the proof.

After the global rescale and prior to the adjustment, we can assume the function $f(r)$ takes the following form in a neighbourhood of $r=a:$
$$f(r)=\left\{ \matrix{(r-a)\cos(R/N)+N\sin(R/N) & r<a\cr\cr
N\sin\Bigl({{R-a+r} \over N} \Bigr) & r \ge a. \cr} \right.$$
The join at $r=a$ is $C^1$ but not smooth. 

Notice that $(\ast)$ holds with equality at all points on the sine curve. For $r<a,$ we have by Theorem 3.2 (2) that $f$ and $f'$ are increasing.
As equality in $(\ast)$ occurs in the limit as $r \rightarrow a^{+},$ we deduce that $(\ast)$ holds strictly for $r<a.$

Smooth $f(r)$ as follows. First we re-choose the joining point to be $r=a+\epsilon$ for some very small $\epsilon>0.$ Replace the function
$f(r)$ above by the following:
$$f(r)=\left\{ \matrix{(r-a)\cos(R/N)+N\sin(R/N) & r<a\cr\cr
\psi(r) & r \in [a,a+\epsilon) \cr\cr
N\sin\Bigl({{R-(a+\epsilon)+r} \over N} \Bigr) & r \ge a+\epsilon \cr} \right.$$
where $\psi(r)$ is any concave down function which joins smoothly with the sine function at $r=a+\epsilon$ and satisfies
$$0<\psi'(r)<\cos\Bigl({{R-(a+\epsilon)+r} \over N} \Bigr)$$ for all $r \in [a,a+\epsilon).$ Such a $\psi$ clearly exists.
For this re-defined function it is clear that $(\ast)$ holds with strict inequality, and in particular at $r=a,$ the single non-smooth point.
It is now clear that we can smooth at this point whilst maintaining $(\ast),$ as such a smoothing can have negligible effect on the values of
the function, while the first derivatives interpolate between the values on either side of the join.

It is evident from the scalar curvature formula in \S4 that the above adjustment to $f$ also maintains positive scalar curvature. 

Having made the adjustment, for notational simplicity we will re-choose the value of $a$ to be $a+\epsilon.$

The rescaled, adjusted metric $g(f,h_{\infty})=dr^2+h_{\infty}(r)ds^2_{n-1}+f(r)^2ds^2_m$ on $[0,a] \times S^{n-1} \times S^m$  
can be viewed as consisting of a $\lq$cap' $([0,a] \times S^{n-1}\,;\,dr^2 + h_{\infty}^2(r)ds^2_{n-1})$
which is topologically a disc $D^n$, and a $\lq$tube' $([0,a] \times S^m\,;\, dr^2+f^2(r)ds^2_m),$ which share the
parameter $r$. Let us consider the issue of smoothly gluing this object into $M \setminus \phi(D^{m+1}_R(N) \times S^{n-1}(\rho)).$
By considering a slightly smaller
disc $D^{m+1}$ if necessary, we may assume that for some $\epsilon'>0,$ a neighbourhood
of the boundary in $M \setminus \phi(D^{m+1} \times S^{n-1})$ takes the form $([R,R+\epsilon'] \times S^{n-1} \times S^m\,;\,
dt^2+\zeta^2(0)ds^2_{n-1}+N^2\sin^2(t/N)ds^2_m)$. We can then split the gluing problem into two: the problem of gluing
the $\lq$cap' onto its $\lq$continuation' $([-\epsilon, 0] \times S^{n-1}\,;\, dr^2+\zeta^2(0)ds^2_{n-1})$; and the
problem of gluing the $\lq$tube' onto its continuation $([R,R+\epsilon'] \times S^m\,;\, 
dt^2+N^2\sin^2(t/N))ds^2_m)$. Clearly, what we need is the smooth joining of the functions $h_{\infty}(r)$ to the constant function with value $\zeta(0)$, and $f(r)$ to $N\sin(t/N)$, when the $r$ and $t$ parameters are suitably concatenated.

We want to go from this surgery picture to adding the handle $\bar{E}$ to $\bar{M}$ in such a way that the boundary metric agrees with that for the surgery, and of course whole metric must be smooth and have positive scalar
curvature. Topologically, going from $E$ to the handle $\bar{E}$
requires the $\lq$filling-in' of the fibre spheres $S^m$, that is, replacing the $\lq$tube' by a corresponding
$\lq$solid tube'. (See Figure 9.)
Again, the smooth gluing problem can be divided into gluing the $\lq$cap' onto its
continuation $([-\epsilon ,0] \times S^{n-1}\,;\, ds^2+\zeta^2(s)ds^2_{n-1})$ and the $\lq$solid tube' onto
its continuation $([-\epsilon ,0] \times D^{m+1} \,;\, ds^2 + dt^2+N^2\sin^2(t/N)ds^2_m).$ Note that in this
handle picture, the sphere radii in the cap are naturally parametrised by $s$, the parameter running along
the central axis in the solid tube. Thus $s$ will be a shared parameter in the handle metric in the same
way that $r$ was a shared parameter for the surgery metric.
\bigskip

\centerline{
\relabelbox \epsfxsize 3.5 truein
\epsfbox{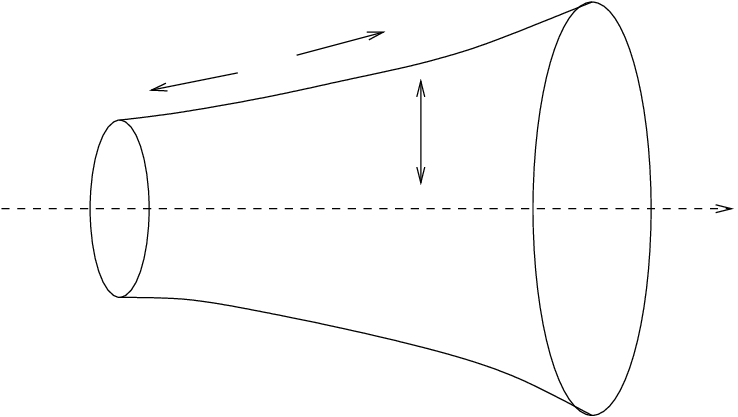} 
\extralabel <-135pt, 80pt> {$f$}
\extralabel <-139pt, 99pt> {$r$}
\extralabel <0pt, 58pt> {$s$}
\endrelabelbox}
\bigskip
\centerline{\it Figure 9: the `solid tube'.}
\bigskip

We focus first on the $\lq$solid tube'. Our approach to defining a metric on this is to view   
it as a subset of a $\lq$big tube' $([0,b] \times D^{m+1}\,;\, ds^2+ dt^2+N^2\sin^2(t/N)ds^2_m)$, for some $b>0$ to be chosen later.
(See Figure 10.) 
\bigskip
\centerline{
\relabelbox \epsfxsize 3.5 truein
\epsfbox{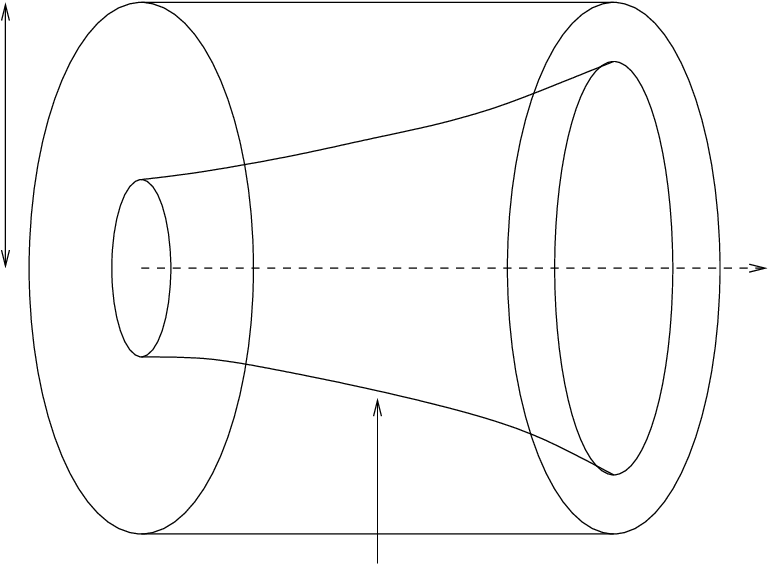} 
\extralabel <-130pt, -10pt> {solid tube}
\extralabel <-220pt, 115pt> {$b$}
\extralabel <-2pt, 79pt> {$s$}
\endrelabelbox}
\bigskip\medskip
\centerline{\it Figure 10: the `big tube'.}
\bigskip
For simplicity, we are using the same symbols $s$ and $t$ for the parameters here as for the corresponding
parameters in $\bar{M},$ on the understanding that a (metrically insignificant) adjustment
will need to be made to the $s$ parameter on one side of the join in order for concatenation to take
place. This $\lq$big tube' will clearly glue onto $([-\epsilon ,0] \times D^{m+1} \,;\, ds^2 + dt^2+N^2\sin^2(t/N)ds^2_m)$ on concatenation, to give a smooth positive scalar curvature metric {\it on the interior}.
We must therefore show how to cut the $\lq$solid tube' from the $\lq$big tube' in such a way that the boundary
of the resulting object is smooth, and has a metric which agrees with that coming from the corresponding surgery.

We will define a curve $\gamma(r)$ in the $(s-t)$-plane (with metric $ds^2+dt^2$) which will define the solid tube. The $s$-coordinate will parameterize the central axis of the solid tube, and the $t$-coordinate will determine
the radius of the boundary spheres $S^m$. We will write $\gamma(r)=(\gamma_s(r),\gamma_t(r))$.

Recall that the (boundary) tube must have induced metric $dr^2 + f^2(r)ds^2_m$. For such a warped product,
the parameter $r$ measures distance {\it around} the tube (as opposed to along the central axis). For this
reason we clearly need $\gamma(r)$ to be a unit speed curve, that is, $\gamma'_s(r)+\gamma'_t(r) \equiv 1$.

We need to choose $\gamma_t(r)$ so that the distance sphere in $(D^{m+1}\,;\,
dt^2+N^2\sin^2(t/N)ds^2_m)$ at a distance of $\gamma_t(r)$ from the centre is round of radius $f(r).$ In other
words we need $N\sin(\gamma_t(r)/N)=f(r),$ and therefore $$\gamma_t(r)=N\sin^{-1} (N^{-1}f(r)).$$
Using this together with $\gamma'_s(r)+\gamma'_t(r) \equiv 1$ allows us to solve for $\gamma_s(r)$: namely
$$\gamma_s(r)=\int_0^r \Bigl( 1-{{f'^2(u)} \over {1-N^{-2}f^2(u)}} \Bigr)^{1 \over 2} \, du.$$ 
Note that the value of $b$ defining the range of $s$ above, is equal to $\gamma_s(a)$, where $r=a$ corresponds
to the end of the tube to be glued onto $\bar{M}.$ (See Figure 11.)
\bigskip

\centerline{
\relabelbox \epsfxsize 2 truein
\epsfbox{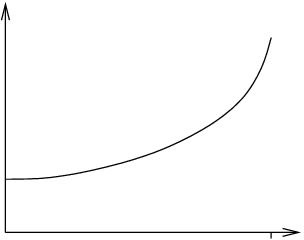} 
\extralabel <-52pt, 57pt> {$\gamma(r)$}
\extralabel <-128pt, 90pt> {$t$}
\extralabel <-2pt, 1pt> {$s$}
\extralabel <-19pt, -9pt> {$b$}
\endrelabelbox}
\bigskip
\centerline{\it Figure 11}
\bigskip\smallskip

Note that the above expression for $\gamma_s(r)$ is well defined because the inequality $(\ast)$ is satisfied, as discussed above.

For notational simplicity, let us denote the solid tube and its metric by $([0,b] \times D^{m+1}\,;\,
ds^2+g_s)$ where $g_s$ is the metric $dt^2+N^2\sin^2(t/N)ds^2_m$ for $t \in [0,\gamma_t(\iota(s))]$, where
$\iota:[0,b] \rightarrow \hbox{\bb R}$ is defined as follows. For $s<b$ the value of $\iota(s)$ is determined by the
equation $\gamma_s(\iota(s))=s.$ Notice that $\iota(s)$ is smooth for these values of $s$ since $\gamma$ is a 
smooth curve. For $s=b$ set $\iota(b)=a,$ which means that $\gamma_t(\iota(b))=R.$ Note that $\iota$ will in
general not be continuous at $s=b.$ The idea here is that the function $\iota(s)$ allows us to replace the parameter $r$
in our subsequent calculations by the parameter $s$ running along the central axis of the solid tube.

We turn our attention now to the $\lq$cap' $D^n$. In the surgery picture, we require the distance sphere at radius
$r$ to have the metric $h_{\infty}^2(r)ds^2_{n-1}.$ In the handle picture, these spheres will be parametrised by $s$, as
discussed above. Clearly, the required metric for any given value of $s$ is $h_{\infty}^2(\iota(s))ds^2_{n-1}.$ Notice that $h(\iota(s))$ is smooth on $[0,b],$ despite the fact that $\iota$ cannot be assumed continuous at $b$. This is because $h_{\infty}(r)$ is constant for $r \in [R',a].$  

A fundamental difference between the handle and surgery cases is that the function 
$h(\iota(s))$ needs to be joined to the function $\zeta(s)$, as opposed to a constant function. However, since $\zeta'>0$ and $\zeta''<0$, this joining can easily be achieved by a small concave down deformation of $h(\iota(s))$ in a neighbourhood of
the joining point. Such a bend can only have a positive effect on the scalar curvature. Note that the second
derivative of $h_{\infty}(\iota(s))$ is non-positive. To see this, recall that $f(r) \equiv 1$ for small $r$,
and this means that for $s$ in the same range we have $\iota(s)=s.$ Also, where $f$ is non-constant we have
$h_{\infty}(r) \equiv c.$ 
The assertion now follows immediately since the desired second derivative is equivalent to $h_{\infty}''(\iota(s))(\iota'(s))^2+h_{\infty}'(\iota(s))\iota''(s).$ 
As a result we deduce that the cap itself has positive scalar curvature (and non-negative Ricci curvature). 

Using the above notation, the whole handle and its metric can now be expressed as $$([0,b] \times S^{n-1} \times D^{m+1}\, , \, ds^2+h_{\infty}^2(\iota(s))+g_s). \eqno{(\ddag)}$$ As the boundary metric by construction agrees with the surgery metric, it is
clear that $\bar{M}\cup \bar{E}$ admits a smooth positive scalar curvature metric, as claimed. 

Finally, we must consider the situation when the connection $\nabla$ is not flat (so the metric $g(f,h_{\infty},\nabla)$
is not a warped product). Clearly, the technique of cutting
out the $\lq$solid tube' from the $\lq$big tube' explained above requires the metric to be a warped product, and so
will not work here. However, recall that the non-flat part of the connection is restricted to the region $r \in (\delta/3,
2\delta/3).$ Accordingly, we are free to assume that $g(f,h_{\infty},\nabla)$ is a warped product
outside this region.

For $r \in [0,\delta]$, recall that the fibres of $E$ are all isometric to a round sphere of radius
$\rho'.$ There is clearly a submersion metric on $\bar{E}$ for $r \in [0, \delta],$ having fibres isometric to $dt^2+N^2 \sin^2(t/N) ds^2_m,$ for which the induced metric on the boundary agrees with that of $E.$
We choose a value for $\rho'$ for which guarantees the Ricci positivity of $E$, and also the Ricci positivity
of $\bar{E}$ restricted to the non-flat region $r \in (\delta/3, 2\delta/3).$ That such a value of $\rho'$
exists follows from Lemma 5.1. 

For $r \ge \delta$ we use the metric constructed for the product
case. This clearly glues smoothly with the above metric for $r \in [0,\delta]$ to create the desired metric on $\bar{E}.$
\hfill{\square}
\medskip

Let $\Sigma$ now be a homotopy sphere representing any element of $bP_{4k}$, and
let $X_p$ be a plumbed manifold as discussed in \S2 with $\Sigma=\partial X_p.$ 
\proclaim Corollary 5.3. The positive scalar curvature metric $g_{p,\infty}$ on $\Sigma$ can be extended across
$X_p$ to give a positive scalar curvature metric $\bar{g}_{p,\infty}$ on $X_p$. \par
\medskip
\noindent {\bf Proof.} The existence of this positive scalar curvature metric on $X_p$ follows immediately from 
repeated application of Proposition 5.2, starting with the Ricci positive metric on the
$\lq$central' disc-bundle guaranteed by Lemma 5.1 and systematically extending the
metric over each plumbing in turn. \hfill{\square}
\medskip
Fixing this metric on $X_p$, consider a neighbourhood of the boundary with normal parameter $w \in [-\epsilon,0]$. By the compactness of $X_p$ there exists an $\epsilon>0$ such that this neighbourhood is diffeomorphic via the exponential map to $[-\epsilon,0] \times \partial X_p.$ We will express the metric here as $dw^2+g(w)$, with $g(w)$
a metric on the equidistant hypersurface at a distance $|w|$ from $\partial X.$ By the openness of the
positive scalar curvature condition, it follows that $\hbox{scal}(g(w))>0$ for all $w$ suitably close
to 0, as well as $\hbox{scal}(dw^2+g(w))>0.$

Fix an arbitrary local coordinate system on $\partial X_p$, $\{x_1,...,x_{4k-1} \}$, and extend to
a local coordinate system on the collar neighbourhood $\{w,x_1,...,x_{4k-1}\}$ via the identification
of the collar with $[-\epsilon,0] \times \partial X_p.$ In the following Lemma, the subscript
$i$ refers to the coordinate system on $\partial X_p$.

\proclaim Lemma 5.4. For $w$ suitably close to 0, the metrics $g(w)$ above satisfy $${{\partial g_{ii}(w)} \over {
\partial w}} \ge 0.$$ \par
\medskip
\noindent {\bf Proof.} To begin with, we will choose the coordinate system $\{x_1,...,x_{4k-1} \}$ to lie entirely on the
boundary of one of the handles $\bar{E}$ used in the construction of $X_p$.

First assume that the metric on $\bar{E}$ takes the form $(\ddag)$, and that the connection is chosen to be flat. We can identify four mutually
orthogonal directions: the $s$-direction with metric $ds^2$; the $t$-direction with metric $dt^2$; the $S^m$ direction
with metric $N^2\sin^2(t/N) ds^2_m;$ and the $S^{n-1}$ direction with metric $h^2_{\infty}(\iota(s)) ds^2_{n-1}.$
The boundary of solid tube is determined by the function $f$, and in turn this determines the direction of the
normal vector $\partial/\partial w.$

Assume the coordinate system lies in the region where $f(\iota(s)) 
\equiv 1.$ Here we have $\partial/\partial w=\partial/\partial t.$ Clearly $g_{ii}$ is non-decreasing along $t$
parameter lines, hence the result.

If the coordinate system lies in the region corresponding where $f(\iota(s))$ is increasing,
we have that $\partial/\partial w$ can be decomposed into non-negative components in the $-\partial/\partial s$
and $\partial/\partial t$ directions. As $h_{\infty}(\iota(s))$ is constant at these values of $s$, we see that
the metric itself is constant along $s$-parameter lines. Again, we have that $g_{ii}$ is non-decreasing along $t$
parameter lines, hence the result in this case.

Now suppose that the connection on $E$ is not everywhere flat. Recall 
that the non-flatness is constrained to lie within the region corresponding to $r \in (\delta/3,2\delta/3)$, which is equivalent to $s \in (\delta/3,2\delta/3).$ Note that changing the principal connection on the associated $SO(m+1)$-bundle only alters the metric on
$\bar{E}$ in directions tangent to the distance sphere bundles. In particular, normal directions to
the distance sphere bundles are unaffected, so in this region $\partial/\partial w=\partial/\partial t$ as before. The same arguments as for the flat case again yield the result that $g_{ii}$ is non-decreasing in normal directions.

Next suppose that our initial coordinate system $\{x_1,...,x_{4k-1} \}$ lies on the disc-bundle
with which we begin our construction of $X_p$, where the metric takes the form described in Lemma 5.1.
Arguments analogous to the above show that $\partial g_{ii}(w)/\partial w \ge 0$,
since $\theta(s)$ is a strictly increasing function.

Finally, if we choose a coordinate system $\{x_1,...,x_{4k-1} \}$ which crosses the boundary
between a handle and the original disc-bundle, or between two handles, the non-negativity of a derivative arbitrarily close to the boundary on either side is clearly sufficient to guarantee non-negativity throughout the region. 
\hfill{\square}
\vfill\eject

\centerline{\bf \S6 Paths of metrics}
\bigskip
In this section we will ultimately prove Theorem A. Firstly, however, we must study the problem of deforming a positive scalar curvature metric on a manifold with boundary to a product in some neighbourhood of the boundary, in such a way that the positive scalar curvature is preserved. For the manifolds and metrics we are concered about in this paper, Lemma 5.4 (above) is crucial to being able to perform such a deformation. The importance of this Lemma arises from the following result:
\proclaim Lemma 6.1. Consider a metric of the form $\eta^2(w)dw^2+g(w)$ on $I \times M$ where $I$ is an interval parametrised by $w$, $M^n$ a manifold, $\eta(w)$ a smooth real-valued function and $g(w)$ a smooth path of metrics on $M$. Given a local coordinate system $\{x_1,...,x_n\}$
for $M$, consider the local coordinate system $\{w,x_1,...,x_n\}$ on $I \times M$. With respect to this coordinate system the scalar curvature of the above metric is given by $$\hbox{scal}(\eta^2(w)dw^2+g(w))=
\hbox{scal}(g(w))+{{H(w)}\over{\eta^2(w)}}+{{\eta'(w)}\over{\eta^3(w)}}\sum_{i=1}^n \Bigl(
{{\partial g_{ii}(w)} \over {\partial w}}/g_{ii}(w) \Bigr)$$ where $H(w)$ is a smooth one-parameter family of real-valued functions on the domain of the coordinate system $\{x_1,...,x_n\}$, which is independent of $\eta(w)$ and given by $$\eqalign{
H(w)=&{1 \over 2}\sum_{i<j} \Bigl( {{\partial g_{ij}} \over {\partial w}} \sum_{k=1}^n \sum_{p=1}^n g^{kp} {{\partial g_{jp}} \over {\partial w}} g_{ki}\,-\, {{\partial g_{jj}} \over {\partial w}} \sum_{k=1}^n \sum_{p=1}^n g^{kp} {{\partial g_{ip}} \over {\partial w}} g_{ki} \Bigr)/(g_{ii}g_{jj}-g^2_{ij}) \cr \cr
&+\sum_{i=1}^n {1 \over {2g_{ii}}}\Bigl( \sum_{k=1}^n \sum_{p=1}^n g^{kp} 
{{\partial g_{ik}} \over {\partial w}} {{\partial g_{ip}} \over {\partial w}}\,-\, 2{{\partial^2 g_{ii}} \over {\partial w^2}} \Bigl).
\cr}$$
\par
\noindent {\bf Proof:} This result follows from an elementary calculation (for example by computing Christoffel symbols for the coordinate system), and we omit the details. \hfill{\square}
\medskip

\proclaim Lemma 6.2. Assume that the metric $dw^2+g(w)$ on $[-\epsilon,0] \times M$ has positive scalar curvature. 
Given a local coordinate system on $M,$ suppose that $$\sum_{i=1}^n {{\partial g_{ii}(w)} \over {\partial w}} \ge 0$$ at all points of $[-\epsilon',0] \times M$ covered by the coordinate system, for some $\epsilon'
\le \epsilon.$ Then there exists a smooth positive function $\beta(w)$ for $w \in [-\epsilon',0]$
with $$\beta(w)=1 \hbox{ for } w \in [-\epsilon',-\epsilon'/2];$$
$$\beta(w)=\Lambda \hbox{ for } w \in [-\epsilon'/4,0]$$ for any given large constant $\Lambda$,
such that the metric $\beta^2(w)dw^2+g(w)$ has positive scalar curvature. \par
\noindent {\bf Proof:} By Lemma 6.1, $\hbox{scal}(dw^2+g(w))=\hbox{scal}(g(w))+H(w)$ for all $w \in
[-\epsilon,0]$, and this is positive by hypothesis. Since $\hbox{scal}(g(w))>0$ it follows that 
$\hbox{scal}(g(w))+H(w)/\beta >0$ for any $\beta \ge 1$. It is then clear that for any {\it increasing}
function $\beta(w)$ of the required form, $$\hbox{scal}(g(w))+{1 \over {\beta^2(w)}}H(w)+
{{\beta'(w)} \over {\beta^3(w)}}\sum_{i=1}^n \Bigl(
{{\partial g_{ii}(w)} \over {\partial w}}/g_{ii}(w) \Bigr)>0$$ since the last term is non-negative by hypothesis.
\hfill{\square}
\medskip

\proclaim Lemma 6.3. If $g(w)$ is any smooth path through positive scalar curvature metrics on a compact manifold $M$, and $w$ belongs to a compact interval $I$, then the metric $\Lambda^2dw^2+g(w)$ has positive scalar curvature on $I \times M$ for all $\Lambda$ sufficiently
large. \par
\noindent {\bf Proof:} By Lemma 6.1, the scalar curvature of $\Lambda^2dw^2+g(w)$ is $\hbox{scal}(g(w))+H(w)/\Lambda^2.$ 
By the compactness of $I \times M$, there exists a constant $C \ge 0$ such that $$\hbox{scal}(\Lambda^2dw^2+g(w)) \ge
\hbox{scal}(g(w))-C/\Lambda^2.$$ Since $\hbox{scal}(g(w))>0$ for all $w \in I$, we can choose $\Lambda$ sufficiently large so that $\hbox{scal}(g(w))-C/\Lambda^2>0,$ and this guarantees the positivity of $\hbox{scal}(\Lambda^2dw^2+g(w)).$ 
\hfill{\square}
\medskip

\proclaim Corollary 6.4. For a given homotopy sphere $\Sigma \in bP_{4k}$, suppose that $X_p$ is a plumbed manifold
as discussed in \S2 with $\partial X_p=\Sigma$. Then $X_p$ can be equipped with a positive scalar curvature metric
which is a product $dw^2+g$ in a neighbourhood of the boundary, and where the component of the moduli space
of positive scalar curvature metrics on $\Sigma$ containing $g$ also contains a Ricci positive metric.
\par
\noindent {\bf Proof:} We equip $X_p$ with the positive scalar curvature metric $\bar{g}_{p,\infty}$ from Corollary 5.3. 
In a neighbourhood of the boundary of $X_p$, this takes the form $dw^2+g(w)$ ($w \le 0$), for some smooth
path of metrics $g(w)$ on $\Sigma$, which for $w$ sufficiently close to zero can be assumed to have positive scalar curvature.
By Lemma 5.4 and Lemma 6.2, we can adjust this metric in a neighbourhood of the boundary to take the form $\Lambda^2dw^2+g(w)$ for any
large constant $\Lambda$, whilst preserving positive scalar curvature. Next, adjoin a collar $\Sigma \times [0,1]$ to
$\partial X_p$ and choose a smooth path of positive scalar curvature metrics $h(w)$ ($w \in [0,1]$) on $\Sigma$ which
smoothly extends the path $g(w)$, and for which $h(w)=h(1)$ for all $w$ sufficiently close to 1. By Lemma 6.3 it follows that
provided $\Lambda$ is chosen sufficiently large, $\Lambda^2dw^2+h(w)$ will have positive scalar curvature. Thus we have
a metric on $X_p \cup (\Sigma \times [0,1])$ which has positive scalar curvature and takes the form $\Lambda^2dw^2+h(1)$ in a
neighbourhood of the boundary. Clearly, a similar statement can be made for $X$ itself, after pulling back this metric
via a suitable diffeomorphism $X_p \cong X_p \cup (\Sigma \times [0,1]),$ and after a global rescale we can assume a positive
scalar curvature metric on $X_p$ taking the required form near the boundary.

By Corollary 4.8, it is clear that $h(1)$ belongs to the same component of positive scalar curvature metrics as the 
Ricci positive metric $g_p$ on the boundary of $X_p$ constructed in [W1] (see Notation 3.5). 
It follows automatically that the same is true when we pass to the moduli space. \hfill{\square}
\medskip

\noindent {\bf Proof of Theorem A.} Regard a homotopy sphere $\Sigma \in bP_{4k}$ as the boundary of a manifold $X_p$ as
described in \S2. As noted at the end of \S2, $\Sigma$ is a spin manifold with vanishing real Pontrjagin classes
and $X_p$ is a parallelisable manifold. Equip $X_p$ with the positive scalar curvature metric guaranteed by Corollary 6.4.
Denote this metric $\bar{g}_p$, and let $g_p$ be its restriction to $\Sigma=\partial X_p$.
Lemma 2.1 applies to the Riemannian manifolds $(\Sigma,g_p)$ and $(X_p,\bar{g}_p)$, and in this case asserts that
$$\eqalign{s(\Sigma,g_p)=&{1 \over {2^{2k+1}(2^{2k-1}-1)}} \sigma(X_p) \cr\cr
&={{8(p|bP_{4k}|+q)} \over {2^{2k+1}(2^{2k-1}-1)}} \cr}$$
for some integer $q$ depending on $\Sigma$. Thus if $p \neq p'$, $|s(\Sigma,g_p)| \neq |s(\Sigma,g_{p'})|,$
and therefore by [KS; Proposition 2.14] $g_p$ and $g_{p'}$ belong to different components of the 
moduli space of positive scalar curvature metrics
on $\Sigma$. It follows immediately that the moduli space of Ricci positive metrics on $\Sigma$  
must have infinitely many components. \hfill{\square}
\bigskip\bigskip\bigskip

\centerline {\bf References}
\bigskip\bigskip
\item{[A]} T. Aubin, {\it M\'etriques riemanniennes et courbure}, J. Differential Geometry {\bf 4} (1970), 383-424.
\bigskip
\item{[Be]} A.L. Besse, {\it Einstein Manifolds}, Springer-Verlag, Berlin (2002).
\bigskip
\item{[BG]} B. Botvinnik, P. Gilkey, {\it The eta invariant and metrics of positive scalar curvature},
Math. Ann. {\bf 302} (1995), no. 3, 507--517. 
\bigskip
\item{[BGN]} C.P. Boyer, K. Galicki, M. Nakamaye, {\it Sasakian geometry, homotopy spheres and positive Ricci curvature}, Topology {\bf 42}
(2003), 981-1002.
\bigskip
\item{[BHSW]} B. Botvinnik, B. Hanke, T. Schick, M Walsh, {\it Homotopy groups of the moduli space of metrics of positive scalar curvature},
Geometry and Topology {\bf 14} (2010), 2047-2076, but see also arXiv:0907.5188.
\bigskip
\item{[Br]} W. Browder, {\it Surgery on simply connected manifolds}, Springer, Berlin, (1972).
\bigskip
\item{[C]} R. Carr, {\it Construction of manifolds of positive scalar curvature}, Trans. Amer. Math. Soc. {\bf 307} no. 1 (1988), 63-74.
\bigskip
\item{[E]} P. Ehrlich, {\it Metric deformations of curvature}, I. Local convex deformations, Geom. Dedicata {\bf 5} (1976), 1-23.
\bigskip
\item{[H]} N. Hitchin, {\it Harmonic spinors}, Advances in Math. {\bf 14} (1974), 1--55.
\bigskip
\item{[JW]} M. Joachim, D. J. Wraith, {\it Exotic spheres and curvature}, Bull. Amer. Math. Soc. {\bf 45} no. 4 (2008), 595-616.
\bigskip
\item{[KM]} M. Kervaire, J. Milnor, {\it Groups of homotopy spheres}, Ann. of Math. {\bf 77} (1963), 504-537.
\bigskip
\item{[KPT]} V. Kapovitch, A. Petrunin, W. Tuschmann, {\it Nonnegative pinching, moduli spaces and bundles with
infinitely many souls}, J. Diff. Geom. {\bf 71} (2005) no. 3, 365-383. 
\bigskip
\item{[KS]} M. Kreck, S. Stolz, {\it Nonconnected moduli spaces of positive sectional curvature metrics}, J. Am. Math.
Soc. {\bf 6} (1993), 825-850.
\bigskip
\item{[LM]} H.B. Lawson, M.-L. Michelsohn, {\it Spin Geometry}, Princeton Math. Series {\bf 38}, Princeton University Press, (1989).
\bigskip
\item{[R]} A. Ranicki, {\it Algebraic and Geometric Surgery}, Oxford University Press (2002).
\bigskip
\item {[SY]} J.-P. Sha, D.-G. Yang, {\it Positive Ricci Curvature on the
Connected Sums of $S^n \times S^m$}, J. Diff. Geom. {\bf 33} (1990), 127-138.
\bigskip
\item{[W1]} D. Wraith, {\it Exotic spheres with positive Ricci curvature}, J. Differential Geometry {\bf 45} (1997), 638-649.
\bigskip
\item{[W2]} D. Wraith, {\it Surgery on Ricci positive manifolds}, J. reine angew. Math. {\bf 501} (1998), 99-113.
\bigskip
\item{[W3]} D. Wraith, {\it Bundle stabilisation and positive Ricci curvature}, Diff. Geom. Appl. {\bf 25} (2007),
552-560.

\bigskip\bigskip
\centerline{D. J. Wraith,}
\centerline{Department of Mathematics,}
\centerline{N.U.I. Maynooth,}
\centerline{Ireland.}
\smallskip
\centerline{email: David.Wraith@nuim.ie} 

\end